\documentclass[11pt, oneside]{amsart}
\usepackage{geometry}
\usepackage{graphicx}	
\usepackage{amssymb}

\title{Prime Gaps in the Gaussian Integers}
\author{Kyle Bradford, James Taylor and George Volz}

\begin{document}
\maketitle

\section*{Introduction}

In elementary number theory, prime numbers emerge as fundamental entities, possessing a unique and intriguing set of properties.  Let $p_{n}$  denote the $n^{th}$  positive, prime integer.  Traditionally, we label the $n^{th}$ prime gap as $g_{n} = p_{n+1} - p_{n}$.  These gaps provide insights into the distribution and clustering of prime numbers, shedding light on their elusive patterns.  

One well-known topic of interest is the asymptotic lower bound of the prime gaps as $n$  grows large.  A famous example is the twin prime conjecture, which suggests that $g_{n} = 2$  infinitely often.  Less well-known but equally significant is the pursuit of tight asymptotic upper bounds of the largest gap.  Cram\'{e}r's conjecture suggests that $g_{n} = O \left( \left( \log p_{n} \right)^{2}  \right)$, providing one of the tightest known upper bounds.

In this paper, we abstract the concept of prime gaps to the Gaussian integer number field $\mathbb{Z}[i]$  using a boxcar metric.  We employ computational methods to find an asymptotic upper bound for our newly-defined prime gaps.  A preliminary section is devoted to developing the necessary tools to describe these prime gaps, followed by computational results and conclusions.The JavaScript code used for these calculations is included in the appendix.

\section*{Preliminaries}

Recall that the boxcar metric on $\mathbb{Z}[i]$  is defined as $d(a+bi,c+di) = |a-c|+|b-d|$.  A ball of radius $r \in \mathbb{Z}^{+}$  centered at $c \in \mathbb{Z}[i]$  is defined as: $B_{r}(c) =\{ x \in \mathbb{Z}[i] : d(x,c) \leq r \}$.  The Gaussian primes in $\mathbb{Z}[i]$  are related to prime integers.  In fact, integers whose absolute values are positive primes congruent to $3 \bmod 4$ are Gaussian primes.  The remaining Gaussian primes $z=a+bi$  in $\mathbb{Z}[i]$  have $a^{2}+b^{2}$ as a prime integer.  

A gap in one dimension translates to a ball in two dimensions, where the ball's interior contains no Gaussian primes, but its boundary contains at least two.  For any Gaussian prime, its complex conjugate is also a Gaussian prime.  Additionally, any unit multiple of this prime or its complex conjugate is also a Gaussian prime, providing an eightfold symmetry. Thus, we restrict our attention to gaps on the set $A= \{ a+bi \in \mathbb{Z}[i] : a \geq b \geq 0 \}$.  

Our goal is to find an asymptotic upper bound for the maximal radii of these balls as we radially move radially away from the origin.  We start with $r=0$  and increment $r$  until  $B_{r}(0) \cap A$  contains at least one Gaussian prime - this is our first gap.  Next, we increment $n$ and consider all $c \in A$  such that $d(c,0) \leq n$.  For each $c$, we again start with $r=0$  and increase $r$  until $B_{r}(c) \cap A$  contains at least one Gaussian prime.  We then take the supremum of these radii over all $c \in A$ with $d(c,0) \leq n$, and consider the asymptotic behavior of the maximal gap as $n$ grows large.  

\section*{Coding}

Thanks to eightfold symmetry, we know we have a lower bound of the X-axis as well as an upper bound of Y = X. This makes for loops perfect for our research-all we need is a starting point and an endpoint. The code is flexible and can start at any X point and grow from there. Once `aStart` and `aMax` are given, a for loop is created, iterating over all X values from `aStart` to `aMax`, incrementing by 1 each loop. We then move into the Y axis, meaning another for loop is needed-this time from 0 to our X value, henceforth 'a'. This allows us to move outward in the field one section at a time, filling in all the space under Y = X. Because points with a Y value of 0 are deduced differently than random points, we must enter an if statement that checks if Y = 0 (henceforth 'b'), and then pass `a` into a separate function that checks if `a` alone is prime.

This function is called `isZeroPrime` and first checks if `a` is less than or equal to 1, because 0 and 1 are not prime. Then, it enters a for loop that increments a testing value `i` by 1. We then say `m = (4 * i + 3)` and check if `a` is equal to `m`. If it is, then our `a` is prime; if not, it loops until `m` is larger than `a` and then returns false. Now, let's return to our gap function. If `isZeroPrime` returns true, we record the point and say it has a gap of 0. Otherwise, we run our point through the `ladder` function, which I'll explain in a moment.

Now, if our b value is not 0, we just check if a*2 + b*2 is prime. We get a number that still must be greater than 1, so that remains the same, but now we need a quick way to determine if this much larger number is prime. We use a for loop that says as long as i is less than the floor of the square root of our number, i (starting at 2) will increase by 1 each loop. If our number divided by i has no remainder, then it's not prime; otherwise, it's prime. Now that you understand the basics of how we checked for primes, let's move on to the ladder function I mentioned earlier.

The ladder function was a difficult function to set up, and it is entered whenever our (a, b) point is not prime-so, quite often. Stick with me while I explain the largest section of code. We bring in our (a, b) values and define them as the center of our sphere. We then want temporary values for `a` and `b` so we can modify them without losing our center. Now, we enter a for loop with a variable called `radius` (R) that increases after each loop. For each `j < R`, `j` increases by 1 (starting at 0). We then say `atemp = a + R - j` and `btemp = b + j` and check if any of those values are prime, making sure to use `isZeroPrime` when necessary. This process repeats for `atemp = a - j` and `btemp = b + R - j`, as well as for `atemp = a - (R - j)` and `btemp = b - j`, and finally for `atemp = a + j` and `btemp = b - (R - j)`. When a prime point is found, we say that the (a, b) point has a distance of `R` to the nearest prime. This is a shortened explanation; if you wish to read the full function, please refer to appendix pages 6-7.

With the ladder function now explained, you understand how we found our gap values. To gather the large number of points we needed, I ran this code for 43 days and collected just under 500 million points. Next, we needed to identify the largest gaps and the largest prime at each increasing radius from zero. To do this, we had the code read the document that contained all of our points. This took some effort, as our computer initially ran out of memory due to the sheer size of the data-6.17 GB. However, once we got it to read, processing the data was easy. We checked whether our point was a certain distance from zero, then compared the gap to the previously saved largest gap. We did the same for the largest prime. This took about a day to process but significantly reduced our data from 500 million points to 31,567 radii from zero.

Lastly, the code found the Big O for Euler from our recorded wedges data. This did not take nearly as long as the wedge calculations since we already had the code for reading the data. We just needed to calculate Euler as log(largest prime)*2 and then calculate O = largest gap / Euler, giving us values of O that seemed to approach zero.

\begin{figure}[h] \label{fig: one}
\includegraphics{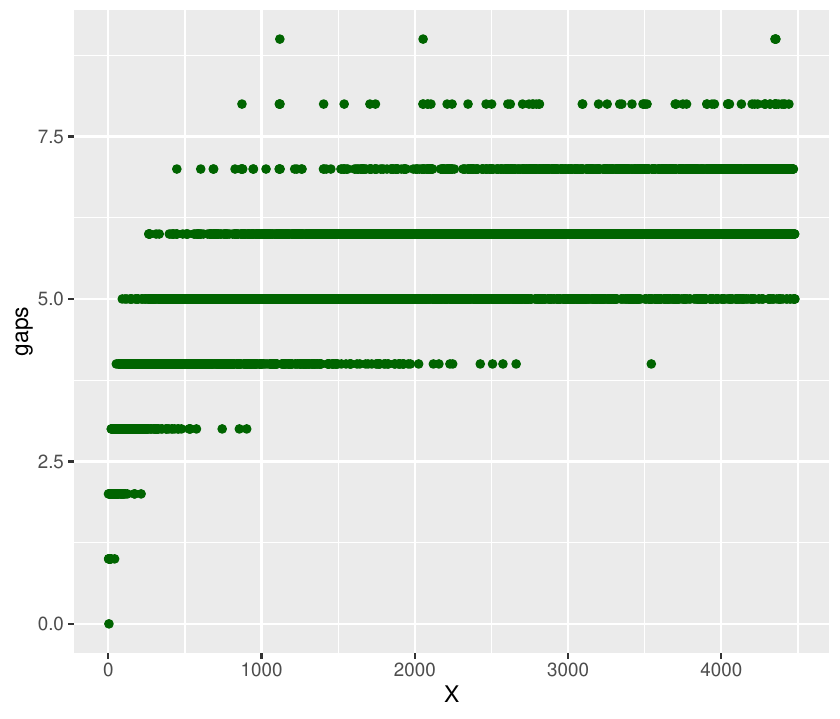}
\caption{X represents the magnitude of the largest prime on the boundary of the gap, and gaps represents the magnitude of the gap}
\end{figure}

\section*{Results}

Our results indicate that, much like the gaps between traditional prime numbers, the maximal gaps between Gaussian primes exhibit a logarithmic pattern as we move radially away from the origin. As you can see in Figure 1 we plot the magnitude of the gaps against the largest magnitude of a prime on the boundary of the gap.  Specifically, we observe a growth rate similar to that suggested by Cram\'{e}r's conjecture for prime gaps in $\mathbb{Z}$, where the gaps increase proportionally to the square of the logarithm of the radius. This suggests that Gaussian primes follow a comparable distributional trend to that of prime numbers in $\mathbb{Z}$, with logarithmic scaling playing a key role in their behavior.

Future work could further refine these bounds and investigate whether more precise analogs to conjectures such as Cram\'{e}r's can be formulated for the Gaussian integers. Nonetheless, our findings confirm that Gaussian primes share fundamental distributional characteristics with their integer counterparts, reinforcing the deep connections between different realms of number theory, while providing a different perspective on what constitutes a gap.

\section*{Appendix}

First, we would like to thank our friend James Taylor for looking over the code and making it more efficient.

\noindent \text{////by: George Volz} \\
\text{//Modified by James Taylor} \\
\text{//This program generates complex numbers and checks if they are Gaussian primes.} \\
\text{//It then finds the gap between each gap.} \\

\noindent \text{import java.util.Scanner;} \\
\text{import java.io.File;} \\
\text{import java.io.FileWriter;} \\
\text{import java.io.IOException;} \\
\text{import java.io.BufferedWriter;} \\
\text{//import java.util.ArrayList;} \\
\text{//import java.util.List;} \\

\noindent \text{public class JavaPrime\{ } \\
\indent \text{ public static void main(String[] args) throws IOException \{ } \\
\indent \qquad \text{Scanner input = new Scanner(System.in);} \\

\indent \qquad \text{System.out.println("1 for yes, 2 for no");} \\
\indent \qquad \text{System.out.print("Are you looking for Primes?");} \\
\indent \qquad \text{int gaps = input.nextInt();} \\
\indent \qquad \text{if(gaps == 1)\{} \\
\indent \qquad \qquad \text{gaps();} \\
\indent \qquad \text{\}} \\
\indent \qquad \text{else\{} \\
\indent \qquad \qquad \text{System.out.print("Are you looking for Wedges?");} \\
\indent \qquad \qquad \text{int wedge = input.nextInt();} \\
\indent \qquad \qquad \text{if(wedge == 1)\{} \\
\indent \qquad \qquad \qquad \text{System.out.println("in the wedge program.");} \\
\indent \qquad \qquad \qquad \text{wedge();} \\
\indent \qquad \qquad \text{\}} \\
            
\indent \qquad \qquad \text{else\{} \\
\indent \qquad \qquad \qquad \text{System.out.print("Are you looking for Euler?");} \\
\indent \qquad \qquad \qquad \text{int euler = input.nextInt();} \\
\indent \qquad \qquad \qquad \text{if(euler ==1)\{} \\
\indent \qquad \qquad \qquad \qquad \text{System.out.println("in the euler program.");} \\
\indent \qquad \qquad \qquad \qquad \text{euler();} \\
\indent \qquad \qquad \qquad \text{\}} \\
\indent \qquad \qquad \text{\}} \\

\indent \qquad \text{\}\}} \\

\indent \text{//gaps functions} \\
\indent \text{public static void gaps() throws IOException\{} \\
\indent \qquad \text{Scanner input = new Scanner(System.in);} \\
\indent \qquad \text{System.out.print("Starting A number: ");} \\
\indent \qquad \text{int aStart = input.nextInt();} \\
\indent \qquad \text{System.out.print("largest number we want to consider: ");} \\
\indent \qquad \text{int aMax = input.nextInt();} \\

\indent \qquad \text{String fileName = "Reaseach Gap Primes.txt";} \\
\indent \qquad \text{int truth = 0;} \\

\indent \qquad \text{for(int a=aStart; a}$<=$ \text{aMax; a++)\{} \\
\indent \qquad \qquad \text{for(int b=0; b}$<=$ \text{a;b++)\{} \\
\indent \qquad \qquad \qquad \text{if (b==0)\{} \\
\indent \qquad \qquad \qquad \qquad \text{if(isZeroPrime(a))\{} \\
\indent \qquad \qquad \qquad \qquad \qquad \text{write(fileName,a,b,truth);} \\
\indent \qquad \qquad \qquad \qquad \qquad \text{continue;} \\
\indent \qquad \qquad \qquad \qquad \text{\}} \\
\indent \qquad \qquad \qquad \qquad \text{else\{} \\
\indent \qquad \qquad \qquad \qquad \qquad \text{ladder(a,b,fileName); // first check arrays]} \\
\indent \qquad \qquad \qquad \qquad \qquad \text{continue;} \\
\indent \qquad \qquad \qquad \qquad \text{\}} \\
                    
\indent \qquad \qquad \qquad \text{\}} \\
\indent \qquad \qquad \qquad \text{if(isPrime(location(a, b)))\{} \\
\indent \qquad \qquad \qquad \qquad \text{write(fileName,a,b,truth);} \\
\indent \qquad \qquad \qquad \text{\}} \\
\indent \qquad \qquad \qquad \text{else\{} \\
\indent \qquad \qquad \qquad \qquad \text{ladder(a,b,fileName); // first check arrays]} \\
\indent \qquad \qquad \qquad \text{\}} \\
\indent \qquad \qquad \text{\}} \\
\indent \qquad \text{\}} \\
\indent \qquad \text{System.out.println("Done!");\}} \\

\indent \text{//checks if a value is prime} \\
\indent \text{public static boolean isPrime(double number) \{} \\
\indent \qquad \text{if (number} $<=$ \text{1) \{} \\
\indent \qquad \qquad \text{return false;} \\
\indent \qquad \text{\}} \\

\indent \qquad \text{for(int i = 2; i} $<=$ \text{(Math.floor(Math.sqrt(number))); i++) \{} \\
\indent \qquad \qquad  \text{if(number\%i==0)} \\
\indent \qquad \qquad \qquad \text{return false;} \\
\indent \qquad \text{\}} \\
\indent \qquad \text{return true;\}} \\

\indent \text{//Check if y=0 x is prime} \\
\indent \text{public static boolean isZeroPrime(double number) \{} \\
\indent \qquad \text{if (number} $<=$ \text{1) \{} \\
\indent \qquad \qquad \text{return false;} \\
\indent \qquad \text{\}} \\

\indent \qquad \text{for(int i = 0; true ; i++) \{} \\
\indent \qquad \qquad \text{int m = (4*i+3);} \\
\indent \qquad \qquad \text{if(number == m)\{} \\
\indent \qquad \qquad \qquad \text{return true;} \\
\indent \qquad \qquad \text{\}} \\
\indent \qquad \qquad \text{if(number}$<$ \text{m)\{} \\
\indent \qquad \qquad \qquad \text{break;} \\
\indent \qquad \qquad \text{\}} \\
\indent \qquad \text{\}} \\
\indent \qquad \text{return false;\}} \\

\indent \text{//writing function} \\
\indent \text{public static void write (String filename, int value1, int value2, int value3)\{} \\
\indent \qquad \text{try \{} \\
\indent \qquad \qquad \text{File f1 = new File(filename);} \\
\indent \qquad \qquad \text{if(!f1.exists()) \{} \\
\indent \qquad \qquad \qquad \text{f1.createNewFile();} \\
\indent \qquad \qquad \text{\}} \\

\indent \qquad \qquad \text{FileWriter fileWritter = new FileWriter(f1.getName(),true);} \\
\indent \qquad \qquad \text{BufferedWriter bw = new BufferedWriter(fileWritter);} \\
\indent \qquad \qquad \text{//for (int i = ((start-minus)+1); i} $<$ \text{(start+1); i++) \{} \\
\indent \qquad \qquad \text{bw.write(value1+"}$|$ \text{"+value2+"}$|$ \text{"+value3);} \\
\indent \qquad \qquad \text{bw.newLine();} \\
\indent \qquad \qquad \text{//\}} \\
\indent \qquad \qquad \text{bw.flush();} \\
\indent \qquad \qquad \text{bw.close();} \\
\indent \qquad \text{\}} \\
\indent \qquad \text{catch(IOException e)\{} \\
\indent \qquad \qquad \text{e.printStackTrace();} \\
\indent \qquad \text{\}\}} \\

\indent \text{//second writing function} \\
\indent \text{public static void writer (String filename, int dist, double largeO)\{} \\
\indent \qquad \text{try \{} \\
\indent \qquad \qquad \text{File f1 = new File(filename);} \\
\indent \qquad \qquad \text{if(!f1.exists()) \{} \\
\indent \qquad \qquad \qquad \text{f1.createNewFile();} \\
\indent \qquad \qquad \text{\}} \\

\indent \qquad \qquad \text{FileWriter fileWritter = new FileWriter(f1.getName(),true);} \\
\indent \qquad \qquad \text{BufferedWriter bw = new BufferedWriter(fileWritter);} \\
\indent \qquad \qquad \text{//for (int i = ((start-minus)+1); i} $<$ \text{(start+1); i++) \{} \\
\indent \qquad \qquad \text{bw.write(dist+"}$|$\text{"+largeO);} \\
\indent \qquad \qquad \text{bw.newLine();} \\
\indent \qquad \qquad \text{//\}} \\
\indent \qquad \qquad \text{bw.flush();} \\
\indent \qquad \qquad \text{bw.close();} \\
\indent \qquad \text{\}} \\
\indent \qquad \text{catch(IOException e)\{} \\
\indent \qquad \qquad \text{e.printStackTrace();} \\
\indent \qquad \text{\}\}} \\

\indent \text{//finding the value of the location} \\
\indent \text{public static double location(int a,int b)\{} \\

\indent \qquad \text{return ((Math.pow(a,2))+(Math.pow(b,2)));\}} \\

\indent \text{//Ladder function} \\
\indent \text{public static void ladder (int centerx, int centery, String FileName)\{} \\
\indent \qquad \text{while(!isPrime(location(centerx, centery)))\{} \\
\indent \qquad \qquad \text{int atemp=0;} \\
\indent \qquad \qquad \text{int btemp=0;} \\
\indent \qquad \qquad \text{for(int radius = 1;true; radius++)\{} \\
\indent \qquad \qquad \qquad \text{//Top right line segment} \\
\indent \qquad \qquad \qquad \text{for (int j=0;j} $<$ \text{radius;j++)\{} \\
\indent \qquad \qquad \qquad \qquad \text{atemp=(centerx+radius-j);} \\
\indent \qquad \qquad \qquad \qquad \text{btemp=(centery+j);} \\
\indent \qquad \qquad \qquad \qquad \text{if(((atemp==0)} $||$ \text{(btemp)==0))\{} \\
\indent \qquad \qquad \qquad \qquad \qquad \text{if((isZeroPrime(atemp)))\{} \\
\indent \qquad \qquad \qquad \qquad \qquad \qquad \text{write(FileName,centerx,centery,radius);} \\
\indent \qquad \qquad \qquad \qquad \qquad \qquad \text{break;} \\
\indent \qquad \qquad \qquad \qquad \qquad \text{\}} \\
\indent \qquad \qquad \qquad \qquad \qquad \text{if((isZeroPrime(btemp)))\{} \\
\indent \qquad \qquad \qquad \qquad \qquad \qquad \text{write(FileName,centerx,centery,radius);} \\
\indent \qquad \qquad \qquad \qquad \qquad \qquad \text{break;} \\
\indent \qquad \qquad \qquad \qquad \qquad \text{\}} \\
\indent \qquad \qquad \qquad \qquad \text{\}} \\
\indent \qquad \qquad \qquad \qquad \text{if(isPrime(location(atemp, btemp)))\{} \\
\indent \qquad \qquad \qquad \qquad \qquad \text{write(FileName,centerx,centery,radius);} \\
\indent \qquad \qquad \qquad \qquad \qquad \text{break;} \\
\indent \qquad \qquad \qquad \qquad \text{\}} \\
\indent \qquad \qquad \qquad \text{\}} \\
\indent \qquad \qquad \qquad \text{if(isPrime(location(atemp, btemp)))} \\
\indent \qquad \qquad \qquad \qquad \text{break;} \\
\indent \qquad \qquad \qquad \text{if(((atemp==0)} $||$ \text{(btemp==0))} \&\& \text{((isZeroPrime(atemp))} $||$ \text{(isZeroPrime(btemp))))} \\
\indent \qquad \qquad \qquad \qquad \text{break;} \\
\indent \qquad \qquad \qquad \text{//Top left line segment} \\
\indent \qquad \qquad \qquad \text{for (int j=0;j}$<$\text{radius;j++)\{} \\
\indent \qquad \qquad \qquad \qquad \text{atemp=(centerx - j);} \\
\indent \qquad \qquad \qquad \qquad \text{btemp=(centery + radius - j);} \\
\indent \qquad \qquad \qquad \qquad \text{if(((atemp==0)}$||$\text{(btemp)==0))\{} \\
\indent \qquad \qquad \qquad \qquad \qquad \text{if((isZeroPrime(atemp)))\{} \\
\indent \qquad \qquad\qquad \qquad \qquad \qquad \text{write(FileName,centerx,centery,radius);} \\
\indent \qquad \qquad \qquad \qquad \qquad \qquad \text{break;} \\
\indent \qquad \qquad \qquad \qquad \qquad \text{\}} \\
\indent \qquad \qquad \qquad \qquad \qquad \text{if((isZeroPrime(btemp)))\{} \\
\indent \qquad \qquad \qquad \qquad \qquad \qquad \text{write(FileName,centerx,centery,radius);} \\
\indent \qquad \qquad \qquad \qquad \qquad \qquad \text{break;} \\
\indent \qquad \qquad \qquad \qquad \qquad \text{\}} \\
\indent \qquad \qquad \qquad \qquad \text{\}} \\
\indent \qquad \qquad \qquad \qquad \text{if(isPrime(location(atemp, btemp)))\{} \\
\indent \qquad \qquad \qquad \qquad \qquad \text{write(FileName,centerx,centery,radius);} \\
\indent \qquad \qquad \qquad \qquad \qquad \text{break;} \\
\indent \qquad \qquad \qquad \qquad \text{\}} \\
\indent \qquad \qquad \qquad \text{\}} \\
\indent \qquad \qquad \qquad \text{if(isPrime(location(atemp, btemp)))} \\
\indent \qquad \qquad \qquad \qquad \text{break;} \\
\indent \qquad \qquad \qquad \text{if(((atemp==0)}$||$\text{(btemp==0))\&\&((isZeroPrime(atemp))}$||$\text{(isZeroPrime(btemp))))} \\
\indent \qquad \qquad \qquad \qquad \text{break;} \\

\indent \qquad \qquad \qquad \text{//Bottom left line segment} \\
\indent \qquad \qquad \qquad \text{for (int j=0;j}$<$\text{radius;j++)\{} \\
\indent \qquad \qquad \qquad \qquad \text{atemp=(centerx-(radius-j));} \\
\indent \qquad \qquad \qquad \qquad \text{btemp=(centery-j);} \\
\indent \qquad \qquad \qquad \qquad \text{if(((atemp==0)}$||$\text{(btemp)==0))\{} \\
\indent \qquad \qquad \qquad \qquad \qquad \text{if((isZeroPrime(atemp)))\{} \\
\indent \qquad \qquad \qquad \qquad \qquad \qquad \text{write(FileName,centerx,centery,radius);} \\
\indent \qquad \qquad \qquad \qquad \qquad \qquad \text{break;} \\
\indent \qquad \qquad \qquad \qquad \qquad \text{\}} \\
\indent \qquad \qquad \qquad \qquad \qquad \text{if((isZeroPrime(btemp)))\{} \\
\indent \qquad \qquad \qquad \qquad \qquad \qquad \text{write(FileName,centerx,centery,radius);} \\
\indent \qquad \qquad \qquad \qquad \qquad \qquad \text{break;} \\
\indent \qquad \qquad \qquad \qquad \qquad \text{\}} \\
\indent \qquad \qquad \qquad \qquad \text{\}} \\
\indent \qquad \qquad \qquad \qquad \text{if(isPrime(location(atemp, btemp)))\{} \\
\indent \qquad \qquad \qquad \qquad \qquad \text{write(FileName,centerx,centery,radius);} \\
\indent \qquad \qquad \qquad \qquad \qquad \text{break;} \\
\indent \qquad \qquad \qquad \qquad \text{\}} \\
\indent \qquad \qquad \qquad \text{\}} \\
\indent \qquad \qquad \qquad \text{if(isPrime(location(atemp, btemp)))} \\
\indent \qquad \qquad \qquad \qquad \text{break;} \\
\indent \qquad \qquad \qquad \text{if(((atemp==0)}$||$\text{(btemp==0))\&\&((isZeroPrime(atemp))}$||$\text{(isZeroPrime(btemp))))} \\
\indent \qquad \qquad \qquad \qquad \text{break;} \\

\indent \qquad \qquad \qquad \text{//bottom right line segment} \\
\indent \qquad \qquad \qquad \text{for (int j=0;j}$<$\text{radius;j++)\{} \\
\indent \qquad \qquad \qquad \qquad \text{atemp=(centerx + j);} \\
\indent \qquad \qquad \qquad \qquad \text{btemp=(centery-(radius-j));} \\
\indent \qquad \qquad \qquad \qquad \text{if(((atemp==0)}$||$\text{(btemp)==0))\{} \\
\indent \qquad \qquad \qquad \qquad \qquad \text{if((isZeroPrime(atemp)))\{} \\
\indent \qquad \qquad \qquad \qquad \qquad \qquad \text{write(FileName,centerx,centery,radius);} \\
\indent \qquad \qquad \qquad \qquad \qquad \qquad \text{break;} \\
\indent \qquad \qquad \qquad \qquad \qquad \text{\}} \\
\indent \qquad \qquad \qquad \qquad \qquad \text{if((isZeroPrime(btemp)))\{} \\
\indent \qquad \qquad \qquad \qquad \qquad \qquad \text{write(FileName,centerx,centery,radius);} \\
\indent \qquad \qquad \qquad \qquad \qquad \qquad \text{break;} \\
\indent \qquad \qquad \qquad \qquad \qquad \text{\}} \\
\indent \qquad \qquad \qquad \qquad \text{\}} \\
\indent \qquad \qquad \qquad \qquad \text{if(isPrime(location(atemp, btemp)))\{} \\
\indent \qquad \qquad \qquad \qquad \qquad \text{write(FileName,centerx,centery,radius);} \\
\indent \qquad \qquad \qquad \qquad \qquad \text{break;} \\
\indent \qquad \qquad \qquad \qquad \text{\}} \\
\indent \qquad \qquad \qquad \text{\}} \\
\indent \qquad \qquad \qquad \text{if(isPrime(location(atemp, btemp)))} \\
\indent \qquad \qquad \qquad \qquad \text{break;} \\
\indent \qquad \qquad \qquad \text{if(((atemp==0)}$||$\text{(btemp==0))\&\&((isZeroPrime(atemp))}$||$\text{(isZeroPrime(btemp))))} \\
\indent \qquad \qquad \qquad \qquad \text{break;} \\
\indent \qquad \qquad \text{\}} \\
\indent \qquad \qquad \text{if(isPrime(location(atemp, btemp)))} \\
\indent \qquad \qquad \qquad \text{break;} \\
\indent \qquad \qquad \text{if(((atemp==0)}$||$\text{(btemp==0))\&\&((isZeroPrime(atemp))}$||$\text{(isZeroPrime(btemp))))} \\
\indent \qquad \qquad \qquad \text{break;} \\
\indent \qquad \text{\}\}} \\

\indent \text{//function for making wedges} \\
\indent \text{public static void wedge() throws IOException\{} \\
\indent \qquad \text{//String reader} \\
\indent \qquad \text{Scanner sc = new Scanner(new File("Reaseach Gap Primes.txt"));} \\
\indent \qquad \text{String name2 = "Wedge Gap research.txt";} \\
\indent \qquad \text{int gaps=0;} \\
\indent \qquad \text{int lprime=0;} \\

\indent \qquad \text{System.out.println("Reading...");} \\
\indent \qquad \text{int num =1;} \\
\indent \qquad \text{int points=1;} \\
\indent \qquad \text{String line = "points";} \\
\indent \qquad \text{//ArrayList}$<$\text{String}$>$ \text{lines = new ArrayList}$<$\text{String}$>$\text{();} \\
\indent \qquad \text{while(sc.hasNextLine())\{} \\
\indent \qquad \qquad \text{line=(sc.nextLine());} \\
\indent \qquad \qquad \text{num++;} \\
\indent \qquad \qquad \text{points++;} \\
\indent \qquad \qquad \text{if(num==10000000)\{} \\
\indent \qquad \qquad \qquad \text{System.out.println("Reading...");} \\
\indent \qquad \qquad \qquad \text{num=1;} \\
\indent \qquad \qquad \text{\}} \\
\indent \qquad \text{\}} \\
\indent \qquad \text{sc.close();} \\
\indent \qquad \text{Scanner scc = new Scanner(new File("Reaseach Gap Primes.txt"));} \\
\indent \qquad \text{//System.out.println(points + " number of points");} \\
\indent \qquad \text{System.out.println("creating the int array");} \\
\indent \qquad \text{int[][] arr = new int[points][3];} \\
\indent \qquad \text{points=0;} \\

\indent \qquad \text{System.out.println("done with preliminary reading.");} \\

\indent \qquad \text{while(scc.hasNextLine())\{} \\
\indent \qquad \qquad \text{String[] temp=scc.nextLine().split("\textbackslash \textbackslash}$|$ \text{");} \\
\indent \qquad \qquad \text{try\{} \\
\indent \qquad \qquad \qquad \text{for(int j = 0; j} $<$ \text{3; j++)\{} \\
\indent \qquad \qquad \qquad \text{arr[points][j] = Integer.parseInt(temp[j]);} \\
\indent \qquad \qquad \qquad \text{\}} \\
\indent \qquad \qquad \qquad \text{points++;} \\
\indent \qquad \qquad \text{\}} \\
\indent \qquad \qquad \text{catch(NumberFormatException ignore)\{\}} \\
     
\indent \qquad \text{\}} \\
\indent \qquad \text{scc.close();} \\
\indent \qquad \text{System.out.println("done reading.");} \\

\indent \qquad \text{for(int d=0;d} $<=$ \text{(arr[points-1][0]);d++)\{} \\
\indent \qquad \qquad \text{for(int k = 0; k} $<$ \text{arr.length; k++)\{} \\
\indent \qquad \qquad \qquad \text{if(d==(arr[k][0]+arr[k][1]))\{} \\
\indent \qquad \qquad \qquad \qquad \text{if(gaps} $<=$ \text{(arr[k][2]))\{} \\
\indent \qquad \qquad \qquad \qquad \qquad \text{gaps=(arr[k][2]);} \\
\indent \qquad \qquad \qquad \qquad \text{\}} \\
\indent \qquad \qquad \qquad \qquad \text{double large = ((Math.pow((arr[k][0]),2))+((Math.pow((arr[k][1]),2))));} \\
\indent \qquad \qquad \qquad \qquad \text{if (arr[k][0]==0)\{} \\
\indent \qquad \qquad \qquad \qquad \qquad \text{if(isZeroPrime(large))\{} \\
\indent \qquad \qquad \qquad \qquad \qquad \qquad \text{if(lprime} $<=$ \text{large)\{} \\
\indent \qquad \qquad \qquad \qquad \qquad \qquad \qquad \text{lprime= (int)large;} \\
\indent \qquad \qquad \qquad \qquad \qquad \qquad \text{\}} \\
\indent \qquad \qquad \qquad \qquad \qquad \text{\}} \\
\indent \qquad \qquad \qquad \qquad \text{\}} \\
\indent \qquad \qquad \qquad \qquad \text{if(isPrime(large))\{} \\
\indent \qquad \qquad \qquad \qquad \qquad \text{if(lprime}$<=$\text{large)\{} \\
\indent \qquad \qquad \qquad \qquad \qquad \qquad \text{lprime= (int)large;} \\
\indent \qquad \qquad \qquad \qquad \qquad \text{\}} \\
\indent \qquad \qquad \qquad \qquad \text{\}} \\
\indent \qquad \qquad \qquad \text{\}} \\
\indent \qquad \qquad \text{\}} \\
\indent \qquad \qquad \text{write(name2,d,gaps,lprime);} \\
\indent \qquad \qquad \text{gaps = 0;} \\
\indent \qquad \qquad \text{lprime = 0;} \\

\indent \qquad \text{\}} \\
        
\indent \qquad \text{System.out.println("Done!");\}} \\
    
\indent \text{//function for making euler} \\
\indent \text{public static void euler() throws IOException\{} \\
\indent \qquad \text{//string scanner} \\
\indent \qquad \text{Scanner sc = new Scanner(new File("Wedge Gap research.txt"));} \\
\indent \qquad \text{String name3 = "Gap V LP Research.txt";} \\
\indent \qquad \text{double pprime=0;} \\
\indent \qquad \text{double euler=0;} \\
\indent \qquad \text{double bigo=0;} \\

\indent \qquad \text{System.out.println("Reading...");} \\
\indent \qquad \text{int num =1;} \\
\indent \qquad \text{int distance=1;} \\
\indent \qquad \text{String line = "distance from Zero";} \\
\indent \qquad \text{//ArrayList}$<$\text{String}$>$ \text{lines = new ArrayList}$<$\text{String}$>$\text{();} \\
\indent \qquad \text{while(sc.hasNextLine())\{} \\
\indent \qquad \qquad \text{line=(sc.nextLine());} \\
\indent \qquad \qquad \text{num++;} \\
\indent \qquad \qquad \text{distance++;} \\
\indent \qquad \qquad \text{if(num==10000000)\{} \\
\indent \qquad \qquad \qquad \text{System.out.println("Reading...");} \\
\indent \qquad \qquad \qquad \text{num=1;} \\
\indent \qquad \qquad \text{\}} \\
\indent \qquad \text{\}} \\
\indent \qquad \text{sc.close();} \\
\indent \qquad \text{Scanner scc = new Scanner(new File("Wedge Gap research.txt"));} \\
\indent \qquad \text{//creating new int array} \\
\indent \qquad \text{int[][] arr = new int[distance][3];} \\
\indent \qquad \text{distance=0;} \\

\indent \qquad \text{System.out.println("done with preliminary reading.");} \\

\indent \qquad \text{while(scc.hasNextLine())\{} \\
\indent \qquad \qquad \text{String[] temp=scc.nextLine().split("\textbackslash \textbackslash}$|$ \text{");} \\
\indent \qquad \qquad \text{try\{} \\
\indent \qquad \qquad \qquad \text{for(int j = 0; j} $<$ \text{3; j++)\{} \\
\indent \qquad \qquad \qquad \text{arr[points][j] = Integer.parseInt(temp[j]);} \\
\indent \qquad \qquad \qquad \text{\}} \\
\indent \qquad \qquad \qquad \text{points++;} \\
\indent \qquad \qquad \text{\}} \\
\indent \qquad \qquad \text{catch(NumberFormatException ignore)\{\}} \\
     
\indent \qquad \text{\}} \\
\indent \qquad \text{scc.close();} \\
\indent \qquad \text{System.out.println("done reading.");} \\

\indent \qquad \text{for(int k = 2; k}$<$\text{(arr.length-1); k++)\{} \\
\indent \qquad \qquad \text{for(int j=0;j}$<=$\text{k;j++)\{} \\
\indent \qquad \qquad \qquad \text{pprime+=(arr[j][2]);} \\
\indent \qquad \qquad \qquad \text{//System.out.println("pprime = " + pprime);} \\
\indent \qquad \qquad \text{\}} \\
\indent \qquad \qquad \text{euler=(Math.pow((Math.log(pprime)),2));}
\indent \qquad \qquad \text{//System.out.println("euler = " + euler);} \\
\indent \qquad \qquad \text{//System.out.println("arr["+k+"][1] = " +arr[k][1]);} \\
\indent \qquad \qquad \text{bigo = ((arr[k][1])/(euler));} \\
\indent \qquad \qquad \text{//System.out.println("tempo is " +tempo);} \\
                
\indent \qquad \qquad \text{writer(name3,k,bigo);} \\
\indent \qquad \qquad \text{bigo = 0;} \\
\indent \qquad \qquad \text{pprime=0;} \\
\indent \qquad \qquad \text{euler=0;} \\
\indent \qquad \text{\}} \\
\indent \qquad \text{System.out.println("Done!");} \\
\indent \text{\}} \\
\noindent \text{\}}

\bibliographystyle{amsplain}

\end{document}